\documentclass[11pt,a4paper,leqno,twoside, headinclude]{amsart}

\title{Cohomological consequences of (almost) free torus actions}
\def\titl{Cohomological consequences of (almost) free torus actions}
\def\auth{Manuel Amann}
\date{April 26th, 2012}

\thanks{The author was supported by a Research Grant of the German Research Foundation.
}

\author{\auth}

\usepackage[english]{babel}

\usepackage[colorlinks,pdftex, plainpages=false]{hyperref}
\hypersetup{pdftitle=\titl, pdfauthor=\auth, pdftoolbar=false,
plainpages=false, hyperindex=true, pdfdisplaydoctitle=true}

\hypersetup{colorlinks,%
citecolor=black,%
filecolor=black,%
linkcolor=black,%
urlcolor=black,%
pdftex}

\usepackage{color}

\usepackage{amsmath, amssymb, amscd, amsthm}
\usepackage{stmaryrd}
\usepackage{latexsym}
\usepackage[all]{xy}
\usepackage{pb-diagram}
\usepackage{rotating}
\usepackage{multicol}
\usepackage{lscape}
\usepackage{wasysym}
\usepackage{longtable}
\usepackage{enumerate}

\xyoption{all}

\newtheorem{theo}{Theorem}[section]
\newtheorem{main}{Theorem}

\newtheorem*{main*}{Theorem}

\newtheorem*{mainprop*}{Proposition}

\newtheorem{mainconj}{Conjecture}

\newtheorem{defi2}[theo]{Definition}
\newtheorem*{defi2*}{Definition}

\newenvironment{defi*}{\begin{defi2*}\normalfont}{\end{defi2*}}

\newenvironment{defin*}[1]{\begin{defi2*}[#1]\normalfont}{\end{defi2*}}
\newtheorem{rem2}[theo]{Remark}
\newenvironment{rem}{\begin{rem2}\normalfont}{\hfill$\boxbox$\end{rem2}}

\newtheorem{lemma}[theo]{Lemma}
\newtheorem{cor}[theo]{Corollary}
\newtheorem*{cor*}{Corollary}

\newtheorem*{conj*}{Conjecture}
\newtheorem*{theo*}{Theorem}
\newtheorem*{ques*}{Question}
\newtheorem*{mi2}{Main Idea}

\newtheorem{ex2}[theo]{Example}

\newtheorem{exer2}[theo]{Exercise}

\newtheorem{alg2}[theo]{Algorithm}

\newcommand{\qq}{{\mathbb{Q}}}                                     
\newcommand{\s}{{\mathbb{S}}}                                      
\newcommand{\zz}{{\mathbb{Z}}}                                     
\newcommand{\B}{{\mathbf{B}}}                                      
\newcommand{\dif} {{\operatorname{d}}}                             
\newcommand{\In} {{\,\subseteq\,}}                                 
\newcommand{\rk}{{\operatorname{rk\,}}}                            
\newcommand{\co}{\colon\thinspace}                                 

\newcommand{\comment}[1]{}                                         
\newcommand{\hto}[1]{\overset{#1}{\hookrightarrow}}                
\newcommand{\step}[1]{\textbf{Step #1.}}                           
\newcommand{\ack}{\noindent\textbf{Acknowledgements. }}            
\newcommand{\str}{\noindent\textbf{Structure of the article. }}    

\newenvironment{prf}{\begin{proof}[\textsc{Proof}]} {\end{proof}}     

\begin{document}

\maketitle \thispagestyle{empty}


\begin{abstract}
The long-standing Halperin--Carlsson conjecture (basically also known as the toral rank conjecture) states that the sum of all Betti numbers of a well-behaved space $X$ (with cohomology taken with coefficients in the cyclic group $\zz_p$ in characteristic $p>0$ respectively with rational coefficients for characteristic $0$) is at least $2^n$ where $n$ is the rank of an $n$-torus $\zz_p \times \stackrel{(n)}\dots \times \zz_p$ (in characteristic $p>0$) respectively $\s^1 \times \stackrel{(n)}\dots \times \s^1$ (in characteristic $0$) acting freely (characteristic $p>0$) respectively almost freely (characteristic $0$) on $X$. This conjecture was addressed by a multitude of authors in several distinct contexts and special cases.

However, having undergone various reformulations and reproofs over the decades, the best known general lower bound on the Betti numbers remains a very low linear one. This article investigates the conjecture in characteristics $0$ and $2$ which results in improving the lower bound in characteristic $0$.
\end{abstract}


\section*{Introduction}


Spurred by the desire to understand the intrinsic beauty of symmetric objects, it is a classical problem in algebraic and differential topology to understand implications of group actions on (compact closed connected) manifolds. This comprises a variety of different techniques and approaches ranging from establishing fixed-point theorems, vanishing theorems for several invariants like the Euler-characteristic or the $\hat A$-genus to equivariant versions of cohomology and homotopy theories, $K$-theory, index theory, etc., just to mention a few examples. The general theory does produce differently flavoured situations basically depending on whether the group action does or does not possess fixed-points.

In the first case several theorems achieve to connect topological properties of the manifold to properties of the fixed-point set. Maybe the most classical example of this kind states that the Euler characteristic of a manifold my be detected as the Euler characteristic of its fixed-point set under an $\s^1$-action.

In this article we shall be concerned with the case when the group action does not possess any fixed-points---more precisely, we shall consider (almost) free group actions of tori (in characteristic $2$ and $0$). The general feeling that in this case the space upon which the torus acts should be large from a point of view of algebraic topology, the idea being that it should contain at least as many ``holes''---in an appropriate sense---as the torus itself is formalised in the classical and long-standing \emph{Halperin--Carlsson conjecture}. Before we state this conjecture, let us introduce the necessary terminology.

By an \emph{$n$-torus in characteristic $p$} we refer to an $n$-fold Cartesian product
\begin{align*}
T=\zz_p \times \stackrel{(n)}\dots \times \zz_p
\end{align*}
where $\zz_p$ denotes the cyclic group of prime order $p$ if $p>0$. If $p=0$ we consider the ``ordinary torus''
\begin{align*}
T=\s^1 \times \stackrel{(n)}\dots \times \s^1
\end{align*}

We recall that an action of a group on a manifold is called \emph{free} if all its isotropy groups are trivial. The action is commonly referred to as being \emph{almost free} if all its isotropy groups are finite.

Let us now state the
\begin{conj*}[Halperin--Carlsson]
If an $n$-torus (in the respective characteristic) acts freely (in characteristic $p$) respectively almost freely (in characteristic $0$) on a finite-dimensional paracompact Hausdorff space $X$, the \emph{total dimension} of its cohomology, i.e.~the sum over all Betti numbers, satisfies the estimate
\begin{align*}
\sum_i \dim_\Bbbk H^i(X,\Bbbk)\geq 2^n
\end{align*}
where $\Bbbk$ is a field of the respective characteristic, $0$ or $p$.
\end{conj*}
The rational version of this conjecture is originally due to Stephen Halperin, the version for finite tori was stated by Gunnar Carlsson. It will be understood in the following that $X$ is always finite-dimensional paracompact Hausdorff.

\vspace{5mm}

Mainly in the rational context, this conjecture essentially is also known as the
\begin{conj*}[toral rank]
If $X$ is a nilpotent finite CW-complex, then
\begin{align*}
\dim_\qq H^*(X,\qq)\geq 2^{\rk_0(X)}
\end{align*}
\end{conj*}
(See \cite[section 7.3.3, p.~283]{FOT08}.)

Here $\rk_0(X)$ denotes the \emph{rational torus rank} of $X$ defined as follows: The non-negative integer $\rk_0(X)$ is the rank of a largest torus $T^n$ such that $T^n$ acts almost freely on a finite CW-complex $Y$ of the same rational homotopy type as $X$. The advantage of this formulation is that one obtains the equivalence of this conjecture with a purely algebraic version of the toral rank conjecture using Sullivan models and the Borel fibration (see \cite[Conjecture 7.21, p.~284, Proposition 7.18, p.~280]{FOT08}).

Let us briefly comment on the importance of the (characteristic $0$) toral rank conjecture within the realm of Rational Homotopy Theory. Here it combines nicely with another famous conjecture by Halperin which is certainly another central conjecture within Rational Homotopy Theory. It states that the rational Leray--Serre spectral sequence of a fibration of nilpotent spaces
\begin{align*}
F\hto{} E \to B
\end{align*}
degenerates at the $E_2$-term provided the fibre is an $F_0$-space, i.e.~rationally elliptic (see \cite[section 32, p.~435]{FHT01}) with positive Euler characteristic. Equivalently, the total dimension of the rational cohomology of $E$ satisfies
\begin{align*}
\dim H^*(E,\qq)=\dim H^*(B,\qq)\cdot \dim H^*(F,\qq)\geq  \dim H^*(F,\qq)
\end{align*}

Due to the equivalence of the algebraic toral rank conjecture and the toral rank conjecture itself, we observe that from a rational viewpoint this conjecture can be restated as follows: Given a fibration of nilpotent spaces
\begin{align*}
T^n\hto{} X \to X/T^n
\end{align*}
is the total dimension of the cohomology of the total space at least as large as the total dimension of the cohomology of the fibre $\dim H^*(T^n,\qq)=2^n$? Thus both conjectures, the Halperin conjecture and the toral rank conjecture, embed naturally into the broad context of estimating cohomology dimensions within fibrations.

\vspace{5mm}

Over the decades the Halperin--Carlsson/toral rank conjecture was studied from a multitude of different perspectives using various different approaches. Let us briefly give a rudimentary overview over some of them. Amongst early contributions we may find \cite{All72}. In \cite{Car80}, \cite{AB88} and \cite{Han09} the conjecture is proven for (special) $\zz_p$-actions for particular $p$ on products of spheres.

The rational version of the conjecture was established for several further special cases like \emph{pure} (see \cite[p.~435]{FHT01}) rationally elliptic spaces---see \cite{AH78}---which comprises the case of biquotients and homogeneous spaces. It was shown to hold for K\"ahler manifolds (respectively, in larger generality, for Hard-Lefschetz manifolds)---see \cite{AH78}. In \cite{JL04} several new lower bounds on the rational torus rank (and particular classification results) are provided in the case of rational two-stage spaces---generalising proofs for the elliptic pure and the coformal case. In the case of a formal quotient space $X/T$ of the fibration $T\hto{} X\to X/T$, the conjecture is proved under an additional assumption in \cite{Mun99}.

There is a Lie algebra version of the toral rank conjecture and partial results in that direction (see \cite{DS88}), there are related questions concerning Hamiltonian bundles as well as several further generalisations in commutative algebra.

For excellent outlines on the conjecture together with simplified proofs we strongly recommend the textbooks \cite{AP93} (in the general case) and \cite{FOT08} in the rational case. See \cite{Pup09} for various further references and contributions.

\vspace{5mm}

Yet, despite all these results, a general lower bound for the sum of all Betti numbers could not be improved beyond first results. In \cite{AP86} it is proved that
\begin{align}\label{eqn07}
\dim_\qq H^*(X,\qq)\geq
\begin{cases}
2n &\textrm{for } n\geq 2\\
2(n+1) &\textrm{for } n\geq 3
\end{cases}
\end{align}
in the case of characteristic $0$. (Partial and simplified) reproofs of this result using different methods and techniques can be found in \cite{AP93}, \cite{FOT08} and \cite{Pup09} for example. However, to the knowledge of the author, throughout all the years no improvement could be made on this lower bound. Despite several reformulations, also in the case of characteristic $p$, the existing lowers bounds do not exceed this rational one and are even worse, in general---see \cite[Theorem 1.1]{Pup09} for a complete picture.

Recently, it was again noted in \cite{Ust12} that the toral rank conjecture is related to a conjecture by Horrock. (See \cite[p.~281]{AP93} for an earlier statement of this kind.) Based upon this observation the author then intended to improve the bound \eqref{eqn07}. However, as a crucial tool the erroneous lemma \cite[Lemma 1.2]{Ust12} is used. The author is indebted to Stephen Halperin, Volker Puppe and Jim Stasheff for independently pointing that out in personal communication.

\vspace{5mm}

The main result of this article is a linear improvement of the lower bound \eqref{eqn07}. This is expressed in
\begin{main}\label{theoA}
Let $\Bbbk$ be a field of characteristic $0$. If an $n$-torus $T$ acts almost freely on a finite-dimensional paracompact Hausdorff space $X$, then
\begin{align*}
\dim H^*(X,\Bbbk)\geq 2(n+\lfloor n/3\rfloor)
\end{align*}
\end{main}
Clearly, $X$ may be taken to be a finite CW-complex or a compact manifold.

\vspace{5mm}

One of the main problems when tackling the Halperin--Carlsson conjecture from an algebraic point of view is to provide an effective way of ``bookkeeping'' the cohomological data. Usually, this seems to be achieved at least partially via showing certain degeneration properties of suitable spectral sequences and subsequent rank considerations on the remaining sheets (for example see \cite{JL04} or \cite{Mun99}). For our goals we shall propose a first approach via graph theory as a main tool for these purposes.

\vspace{3mm}

\str In section \ref{sec01} we shall provide the main technical tools for proving theorem \ref{theoA} before we set this information into its topological fundament in section \ref{sec02}---thus finishing the proof of the theorem.

\vspace{3mm}

\ack
The author is very grateful to Jim Stasheff for various inspiring discussions; the author is indebted to Volker Puppe and Jim Stasheff for commenting on a previous version of this article.


\section{An improved lower bound}\label{sec01}

This section is devoted to providing the main technical results for proving theorem \ref{theoA}. As long as possible we shall try to treat characteristics $2$ and $0$ in a similar way leaving room for possible further improvements in characteristic $2$.
By $\Bbbk$ we shall denote a field of characteristic $2$ or of characteristic $0$.
Inspired by the excellent work in \cite{Pup09} we shall respect the notation used there:

Let $R:=\Bbbk[t_1,\dots, t_n]$ be the polynomial algebra over $n$ generators and denote by $\Lambda\langle \cdot \rangle$ the exterior algebra (over $\Bbbk$) over the respective specified generators. We shall decorate elements with degrees. As a convention in characteristic $2$ we set $\deg t_i=1$ and $\deg s_i^m=m$. In characteristic $0$ we agree on $\deg t_i=2$ and $\deg s_i^m=2m+1$. Note that superscripts on the $t_i$ indicate powers, superscripts on the $s_i^m$ refer to degree in the depicted way. (This completely mimics the notation used in \cite{Pup09}.)

We form the \emph{Koszul complex} $K_n(m)$ (for $m\geq 0$) belonging to the regular sequence $(t_1^m,\dots, t_n^m)$ in $R$, i.e.~the commutative differential graded algebra
\begin{align*}
K_n(m):=(\Bbbk[t_1,\dots, t_n] \otimes \Lambda\langle s_1^m, \dots, s_n^m\rangle,\dif)
\end{align*}
with $\dif$ defined by $\dif t_i=0$ for $1\leq i\leq n$, by $\dif s_i^m=t_i^{m+1}$ and extended as a derivation to the entire algebra. By the multi-index $s_{i_1,i_2,\dots}^m$ we denote the element $s_{i_1}\cdot s_{i_2}\cdots \in \Lambda\langle s_1^m, \dots, s_n^m\rangle$.

We consider an $R$-linear map $\gamma$ of differential modules
\begin{align*}
\gamma\co K_n(m)\to K_n(0)
\end{align*}
which lifts the projection
\begin{align*}
H(K_n(m))=R/(t_1^{m+1},\dots, t_n^{m+1})\to H(K_n(0))=R/(t_1^1,\dots, t_n^1)\cong \Bbbk
\end{align*}
In particular, this implies that $\gamma$ is unital, i.e.~$\gamma(1)=1$.

By its \emph{rank} $\rk(\gamma)$ we denote the rank of $\gamma \otimes R_{(0)}$ as a map of $R_{(0)}$-vector spaces. Here, by $R_{(0)}$ we denote the field of fractions of $R$, the localisation of $R$ in the variables $t_i$; i.e.~we adjoin formal inverses to the variables $t_i$.

\vspace{5mm}

We shall now carry out an analysis of the map $\gamma$ in order to find a lower bound on its rank. Our main tool for this will be constituted by graph theory.
Thus let us briefly recall the trivial lemma \ref{lemma04} below, which we slightly adapt for our purposes and which will serve as a simple but very useful tool.

A \emph{sink} of a directed graph is a vertex without outgoing arrows. A \emph{directed acyclic graph} is a directed graph without any directed cycles. Let us modify this terminology by calling a directed graph \emph{directed $l$-acyclic} if it has no directed cycles of length larger than or equal to $l$. For example, a $3$-acyclic graph is acyclic up to subgraphs of the form
$
\xygraph{
!{(0,0) }*+{\bullet_{}}="a"
!{(1,0) }*+{\bullet_{}}="b"
"a":@/_/"b"
"b":@/_/"a"
}
$. Hence a directed graph is acyclic if and only if it is $2$-acyclic.

Correspondingly, we define a \emph{$3$-sink} as a vertex which is a sink up to possibly one $2$-cycle, i.e., more precisely, a vertex with only ingoing edges except for at most one double/bidirectional edge.
\begin{lemma}\label{lemma04}
\begin{itemize}
\item
Let $\mathcal{G}$ be a finite directed acyclic graph. Then $\mathcal{G}$ has a sink.
\item
Let $\mathcal{G}$ be a finite directed $3$-acyclic graph. Then $\mathcal{G}$ has a $3$-sink.
\end{itemize}
\end{lemma}
\begin{prf}
Start in one vertex and follow a directed edge to another vertex. Continue in this fashion. The fact that the graph is acyclic guarantees that no vertex is visited twice. Since the graph is finite, one ends in a vertex without outgoing edges---a sink.

The proof of the second part of the lemma proceeds in the analogous way: Follow directed edges unless they lead back to the vertex preceding the current one in the walk, i.e.~ignoring ``return trips''. The $3$-acyclicity guarantees ending in a $3$-sink.
\end{prf}

\vspace{5mm}

Graph theory will provide an effective way of bookkeeping some of the algebraic data needed to describe the image of $\gamma$. Indeed, lemma \ref{lemma04} will allow us to ``localise'' combinatorical problems arising in the following lemma.

The morphism $\gamma$ is not assumed to be degree-preserving in
\begin{lemma}\label{lemma02}
In both characteristic $2$ and $0$, the map $\gamma$ is injective on the $R$-module
\begin{align*}
\langle \dif s_{1,2,3}^m, \dif s_{4,5,6}^m, \dif s_{7,8,9}^m, \dots\rangle_R
\end{align*}
\end{lemma}
\begin{prf}
\step{1} Since $\gamma$ is unital, $R$-linear and commutes with differentials, we obtain that
\begin{align*}
\gamma(t_i)=t_i
\intertext{which implies that}
\gamma(s_i^m)=t_i^m s_i^0+r_i
\end{align*}
where $\dif r_i=0$.

Moreover, note that $\gamma|_{R\otimes \Lambda^1\langle s_1^m,\dots,s_r^m\rangle}$ is injective.
For this we quote the respective observation in the proof of \cite[Lemma 2.1.b]{Pup09}.

\vspace{5mm}

Let us now consider the terms $\gamma(s_{i,j}^m)$. We obtain
\begin{align*}
\gamma\co \dif(s_{i,j}^m)=t_i^{m+1}s_j^m\pm t_j^{m+1}s_i^m\mapsto t_i^{m+1}t_j^m s_j^0+t_i^{m+1}r_j\pm t_j^{m+1} t_i^m s_i^0\pm t_j^{m+1}r_i
\end{align*}
Since $\dif$ reduces the word-length in the $s_i^m$ by exactly $1$, we obtain that $\gamma(s_{i,j}^m)$ must contain summands of word-length $2$ in the $s_i^0$. Thus
\begin{align*}
\gamma(s_{i_0,j_0}^m)=\sum_{i,j} p_{i,j}s_{i,j}^0+q_{i_0,j_0}
\end{align*}
with $p_{i,j}$ polynomials in the $t_i$ and $q_{i_0,j_0}$ an expression only containing summands which consist of terms with word-length other than $2$ in the $s_i^0$. We therefore obtain that
\begin{align}
\nonumber&\gamma(\dif(s_{i_0,j_0}^m))|_{R\otimes \Lambda^1 \langle s_1^0,\dots,s_n^0\rangle}
\nonumber\\=&\dif(\gamma(s_{i_0,j_0}^m))|_{R\otimes \Lambda^1 \langle s_1^0,\dots,s_n^0\rangle}
\nonumber\\=&\dif\bigg(\sum_{i,j} p_{i,j}s_{i,j}^0\bigg)
\nonumber\\=&\sum_{i,j} p_{i,j} (t_is_j^0\pm t_js_i^0)\nonumber
\intertext{whilst}
&\gamma(\dif(s_{i_0,j_0}^m))|_{R\otimes \Lambda^1 \langle s_1^0,\dots,s_r^0\rangle}
\nonumber\\=&\gamma(t_{i_0}^{m+1}s_{j_0}^m\pm t_{j_0}^{m+1}s_{i_0}^m)|_{R\otimes \Lambda^1 \langle s_1^0,\dots,s_r^0\rangle}
\nonumber\\=&(t_{i_0}^{m+1}t_{j_0}^ms_{j_0}^0+t_{i_0}^{m+1}r_{j_0}\pm t_{j_0}^{m+1}t_{i_0}^ms_{i_0}^0\pm t_{j_0}^{m+1}r_{i_0})|_{R\otimes \Lambda^1 \langle s_1^0,\dots,s_r^0\rangle}
\nonumber\\=&t_{i_0}^{m}t_{j_0}^m(t_{i_0}s_{j_0}^0\pm t_{j_0}s_{i_0}^0)+(t_{i_0}^{m+1}r_{j_0}\pm t_{j_0}^{m+1}r_{i_0})|_{R\otimes \Lambda^1 \langle s_1^0,\dots,s_r^0\rangle}\label{eqn01}
\end{align}
where $|_{R\otimes \Lambda^1 \langle s_1^0,\dots,s_r^0\rangle}$ denotes the projection onto the subspace \linebreak[4] $R\otimes \Lambda^1 \langle s_1^0,\dots,s_r^0\rangle$. Comparing the coefficients in $s_{i_0}^0$ and $s_{j_0}^0$ we see that we actually obtain
\begin{align}\label{eqn02}
\gamma(s_{i_0,j_0}^m)=\sum_{i,j} p_{i,j}s_{i,j}^0+q_{i_0,j_0}
\end{align}
with $p_{i_0,j_0}\neq 0$. Indeed, we realise that every term $p_{i,j}s_{i,j}^0$ under $\dif$ yields a coefficient $t_i$ for $s_j^0$. Thus there has to be a non-trivial term of the form $p_{i_0,j_0}s_{i_0,j_0}^0=t_{i_0}^mt_{j_0}^ms_{i_0,j_0}^0$ in \eqref{eqn02}, since only the terms $t_{i_0}^m$ and $t_{j_0}^m$ arise as coefficients of $s_{j_0}^0$ in \eqref{eqn01}---clearly, terms of the form $s_{i,i}^0$ vanish in the exterior algebra formed by the $s_i^0$.

\vspace{5mm}

Using this insight we compute
\begin{align*}
&\gamma(\dif s_{i_1,i_2,i_3}^m)|_{R\otimes \Lambda^1 \langle s_1^0,\dots,s_n^0\rangle}\\
=&\gamma(t_{i_1}^{m+1} s_{i_2,i_3}^m\pm t_{i_2}^{m+1} s_{i_1,i_3}^m\pm t_{i_3}^{m+1} s_{i_1,i_2}^m)|_{R\otimes \Lambda^2 \langle s_1^0,\dots,s_n^0\rangle}
\\=&t_{i_1}^{m+1}t_{i_2}^mt_{i_3}^m s_{i_2,i_3}^0\pm t_{i_2}^{m+1}t_{i_1}^mt_{i_3}^m s_{i_1,i_3}^0 \pm t_{i_3}^{m+1}t_{i_1}^mt_{i_2}^m s_{i_1,i_2}^0 + t_{i_1}^{m+1} r_{i_2,i_3} \\& + t_{i_2}^{m+1}r_{i_1,i_3}+t_{i_3}^{m+1} r_{i_1,i_2}
\end{align*}
where the $r_{i,j}\in R\otimes \Lambda^2 \langle s_1^0,\dots,s_n^0\rangle$ have word-length $2$ in the $s_i^0$. Note that $\dif(t_{i_1}^{m+1} r_{i_2,i_3}+ t_{i_2}^{m+1}r_{i_1,i_3}+t_{i_3}^{m+1} r_{i_1,i_2})=0$, since $\gamma$ commutes with differentials.

\vspace{5mm}

\step{2} We shall now show that an arbitrary non-trivial linear combination
\begin{align}\label{eqn05}
v=\sum_{i} p_{i} \dif s_{i_1,i_2,i_3}^m
\end{align}
with $p_{i}\in R$ of the specified $\dif s_{i_1,i_2,i_3}$, i.e.~$(i_1,i_2,i_3)\in \{(1,2,3),(4,5,6),\dots\}$, has non-trivial image under $\gamma$.

We compute that
\begin{align*}
\gamma(v)|_{R\otimes \Lambda^2\langle s_1^0,\dots,s_n^0\rangle}=&\sum_i p_i (t_{i_1}^{m+1} (t_{i_2}^mt_{i_3}^m s^0_{i_2,i_3}+r_{i_2,i_3})\\&\pm t_{i_2}^{m+1}(t_{i_1}^mt_{i_3}^m s^0_{i_1,i_3}+r_{i_1,i_3})\pm t_{i_3}^{m+1}(t_{i_1}^mt_{i_2}^m s^0_{i_1,i_2}+r_{i_1,i_2}))
\end{align*}
with
$t_{i_1}^{m+1}r_{i_2,i_3}\pm t_{i_2}^{m+1}r_{i_1,i_3}\pm t_{i_3}^{m+1}r_{i_1,i_3}\in \ker \dif\cap R\otimes \Lambda^2\langle s_1^0,\dots,s_n^0\rangle$.
We shall call the terms $t_{i_1}^{m+1}t_{i_2}^mt_{i_3}^m s_{i_2,i_3}^0$ ``regular terms'' and the terms $t_{i_1}^{m+1}r_{i_2,i_3}$ ``rest terms'' in the following.

\vspace{5mm}

In order to prove that $\gamma(v)|_{R\otimes \Lambda^2\langle s_1^0,\dots,s_n^0\rangle}\neq 0$ we assume the contrary and lead it to a contradiction. So $\gamma(v)=0$ and each of its summands is cancelled by another one. Since the $s_{i_1,i_2,i_3}^0$ are formed over pairwise disjoint triples $(i_1,i_2,i_3)$ it is obvious that a ``regular term'' can only be cancelled by (a summand of a) ``rest term'' (even up to multiples in $R$).

Moreover, a ``regular term'' cannot be cancelled by a ``rest term'' belonging to the same triple $(i_1,i_2,i_3)$. This follows directly from degree reasons and the observations made with equation \eqref{eqn02}. More precisely, if the term $t_{i_1}^{m+1}t_{i_2}^mt_{i_3}^m s^0_{i_2,i_3}$ is cancelled by one of the terms belonging to the triple $(i_1,i_2,i_3)$, by degree reasons in $t_i^{m+1}$ this can only be the ``rest term'' $t_{i_1}^{m+1}r_{i_2,i_3}$. However, such a cancellation would imply that \linebreak[4]$\gamma(s_{i_2,i_3}^m)=\sum_{i,j} p_{i,j}s_{i,j}^0+q_{i_0,j_0}$ with $p_{i_2,i_3}=0$ contradicting \eqref{eqn02}.

\vspace{5mm}

Let us now consider how different ``regular terms'' may be cancelled by ``rest terms'' in general. Since several identical summands (even up to multiples in $R$) might arise, it is clearly not uniquely determined which summands are considered to  cancel each other. However, we choose once and for all an arbitrary but then fixed ``cancellation scheme''.

Using this, we form a directed graph $\mathcal{G}=(V,E)$ with vertices $V=\{x_i\}$ corresponding to the $s_{i_1,i_2,i_3}^m$ with $(i_1,i_2,i_3)\in \{(1,2,3),(4,5,6),(7,8,9),\dots\}$.
We now set a directed edge $(x_i,x_j)=(s_{i_1,i_2,i_3}^m,s_{j_1,j_2,j_3}^m)\in E$ (expressed by the \emph{ordered} pair) between two vertices $x_i=s_{i_1,i_2,i_3}^m$ and $x_j=s_{j_1,j_2,j_3}^m$ if there is a cancellation between a ``regular term'' of $\gamma(\dif s_{i_1,i_2,i_3}^m)$ with a (summand of a) ``rest term'' of $\gamma(\dif s_{j_1,j_2,j_3}^m)$, i.e.~the existence of a directed edge between $s_{i_1,i_2,i_3}^m$ and $s_{j_1,j_2,j_3}^m$ tells us exactly that one of the summands $t_{i_1}^{m+1}t_{i_2}^mt_{i_3}^m s^0_{i_2,i_3}$, $t_{i_2}^{m+1}t_{i_1}^mt_{i_3}^m s^0_{i_1,i_3}$ or $t_{i_3}^{m+1}t_{i_1}^mt_{i_2}^m s^0_{i_1,i_2}$ is cancelled by a (summand of a) ``rest term'' $t_{j_1}^{m+1}r_{j_2,j_3}$, $t_{j_2}^{m+1}r_{j_1,j_3}$ or $t_{j_3}^{m+1}r_{j_1,j_2}$ (always understood up to multiples in $R$). As we have seen, this graph completely encodes how ``regular terms'' may cancel. In particular, it does not possess any loops.

It is our goal to show that not all ``regular terms'' of $\gamma(v)|_{R\otimes \Lambda^2\langle s_1^0,\dots,s_n^0\rangle}$ can be cancelled. This will produce a contradiction proving the injectivity of $\gamma$. We shall do so by establishing two statements which provide a way to treat this problem ``locally''.
\begin{enumerate}
\item \label{case2}
\begin{enumerate}
\item \label{case2-1}
In an induced subgraph
$
\xygraph{
!{<0cm,0cm>;<1.5cm,0cm>:<0cm,1.2cm>::}
!{(0,0) }*+{\bullet_{x_1}}="a"
!{(1,0) }*+{\bullet_{x_2}}="b"
"a":"b"
}
$
of $\mathcal{G}$ consisting of one edge joining two vertices at most $3$ (out of $6$) ``regular terms'' can cancel; namely the ones corresponding to $x_1$.
\item \label{case2-2}
In an induced subgraph
$
\xygraph{
!{(0,0) }*+{\bullet_{x_1}}="a"
!{(1,0) }*+{\bullet_{x_2}}="b"
"a":@/_/"b"
"b":@/_/"a"
}
$
of $\mathcal{G}$ consisting of two edges joining two vertices at most $2$ (out of $6$) ``regular terms'' can cancel, one for each vertex.
\end{enumerate}
In neither case all ``regular terms'' of $x_2$ can be cancelled out.
\item \label{case1} The graph $\mathcal{G}$ is a directed $3$-acyclic graph.
\end{enumerate}

Given these two assertions the injectivity of $\gamma$ on the specified subspace follows: Since $\mathcal{G}$ is $3$-acyclic, by lemma \ref{lemma04} it follows that it possesses a $3$-sink $x$. According to the differentiation provided in assertion \eqref{case2} at most one out of the three ``regular terms'' corresponding to the vertex $x$ may be cancelled (using coefficients from $R$). Thus, in total, $\gamma(v)\neq 0$; a contradiction.

\vspace{5mm}

\step{3} It remains to prove the two auxiliary assertions.

\textbf{ad \eqref{case2}} Let us investigate how many terms may cancel in a subgraph of $\mathcal{G}$ induced by two vertices $x_1$ and $x_2$. Denote by
\begin{align*}
&t_{i_1}^{m+1} (t_{i_2}^mt_{i_3}^m s^0_{i_2,i_3}+r_{i_2,i_3})\pm t_{i_2}^{m+1}(t_{i_1}^mt_{i_3}^m s^0_{i_1,i_3}+r_{i_1,i_3})\\&\pm t_{i_3}^{m+1}(t_{i_1}^mt_{i_2}^m s^0_{i_1,i_2}+r_{i_1,i_2})
\intertext{the terms corresponding to $x_1$ and by}
&t_{j_1}^{m+1} (t_{j_2}^mt_{j_3}^m s^0_{j_2,j_3}+r_{j_2,j_3})\pm t_{j_2}^{m+1}(t_{j_1}^mt_{j_3}^m s^0_{j_1,j_3}+r_{j_1,j_3})\\&\pm t_{j_3}^{m+1}(t_{j_1}^mt_{j_2}^m s^0_{j_1,j_2}+r_{j_1,j_2})
\end{align*}
the ones of $x_2$. By assumption, the sets $\{i_1,i_2,i_3\}$ and $\{j_1,j_2,j_3\}$ are disjoint. This implies that for any terms to cancel, without restriction, the terms corresponding to $x_1$ have to be multiplied by a multiple of $t_{j_1}^{m+1}$. However---we cancel out common multiples for $x_1$ and $x_2$---this implies that the terms $\pm t_{j_2}^{m+1}t_{j_1}^mt_{j_3}^m s^0_{j_1,j_3}$ and $\pm t_{j_3}^{m+1}t_{j_1}^mt_{j_2}^m s^0_{j_1,j_2}$ remain---irrespective of chosen coefficients from $R$.

The two subcases of \eqref{case2} now arise from the following distinction:
For the term $t_{j_1}^{m+1} t_{j_2}^mt_{j_3}^m s^0_{j_2,j_3}$ to be cancelled the coefficient $p_j$ must be divisible by one of $t_{i_1}^{m+1}, t_{i_2}^{m+1}, t_{i_3}^{m+1}$, say $t_{i_1}^{m+1}$ without restriction. If this is the case, the ``regular terms'' $\pm t_{i_2}^{m+1}t_{i_1}^mt_{i_3}^m s^0_{i_1,i_3}$ and $\pm t_{i_3}^{m+1}t_{i_1}^mt_{i_2}^m s^0_{i_1,i_2}$ do not cancel. Thus either at most all of the ``regular terms'' of $x_1$ and none of the ``regular terms'' of $x_2$ cancel---yielding case \eqref{case2-1}---or at most one ``regular term'' of $x_1$ and at most one of $x_2$ cancel---corresponding to case \eqref{case2-2}.

\vspace{5mm}

\textbf{ad \eqref{case1}.}
Let us assume that $\mathcal{G}$ possesses a subgraph of the form
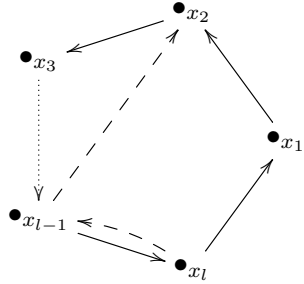
\begin{align*}
\xygraph{
!{<0cm,0cm>;<1cm,0cm>:<0cm,1cm>::}
!{(0,0);a(0)**{}?(1.8)}*+{\bullet_{x_1}}="b1"
!{(0,0);a(72)**{}?(1.8)}*+{\bullet_{x_2}}="b2"
!{(0,0);a(144)**{}?(1.8)}*+{\bullet_{x_3}}="b3"
!{(0,0);a(216)**{}?(1.8)}*+{\bullet_{x_{l-1}}}="b4"
!{(0,0);a(288)**{}?(1.8)}*+{\bullet_{x_l}}="b5"
"b1":"b2" "b2":"b3"  "b4":"b5" "b5":"b1" "b5":@{-->}@/_/"b4" "b4":@{-->}"b2"
"b3":@{..>}"b4"
}
\end{align*}
i.e.~a directed cycle of lengths $l\geq 3$. (This subgraph is not necessarily an induced subgraph, i.e.~we ignore bidirected edges and possibly further edges between the $x_i$ in this schematic---we try to indicate that by dashed arrows, whereas the dotted one stands for the remaining elements of the cycle.) In other words, a ``regular term'' of each $x_{k-1}$ cancels with a (summand of a) ``rest term'' of $x_{k}$ (with suitably chosen $p_k$, $p_{k-1}$) for $1\leq k\leq l$ and the same holds true for $x_l$ together with $x_1$.

The case by case check for the proof of assertion \eqref{case2} yields inductively---irrespective of whether bidirectional edges occur or not, i.e.~in both cases \eqref{case2-1} and \eqref{case2-2}---the following: For the described cancellation to be realised the terms of $x_1$ have to be multiplied by $t_{j''}^{m+1}$ where $j''$ is one of the integers belonging to $x_2$. The terms of $x_2$ are again multiplied by a $t_{j'''}^{m+1}$ with $j'''$ corresponding to $x_3$, etc. However, then for the cancellation between $x_1$ and $x_2$ to be maintained, also the terms corresponding to $x_1$ have to be multiplied with $t_{j'''}^{m+1}$---here again we make use of the fact that the integral triples belonging to the different vertices are disjoint as sets. That is, multiplication with these $t_j^{m+1}$-factors is transitive through the directed graph. The fact that there is an edge
$
\xygraph{
!{<0cm,0cm>;<1.5cm,0cm>:<0cm,1.2cm>::}
!{(0,0) }*+{\bullet_{x_l}}="a"
!{(1,0) }*+{\bullet_{x_1}}="b"
"a":"b"
}
$ then shows iteratively via this transitivity that the terms corresponding to $x_1$ have to be a $t_{j'}^{m+1}$-multiple of themselves, where $j'$ is an integer belonging to $x_1$. This is clearly not possible. Thus a closed unidirectional path of length at least $3$ cannot exist and the graph $\mathcal{G}$ is directed $3$-acyclic. This finishes the proof of the lemma.

\end{prf}

Although we shall not make use of this, we shall state the following observation: The methods of the proof of lemma \ref{lemma02} directly generalise to give the following extension of it.
\begin{lemma}\label{lemma06}
In both characteristic $2$ and $0$ the map $\gamma$ is injective on the $R$-modules
\begin{align*}
\langle \dif s_{1,2,3,4}^m, \dif s_{5,6,7,8}^m, \dif s_{9,10,11,12}^m, \dots\rangle_R\\
\langle \dif s_{1,2,3,4,5}^m, \dif s_{6,7,8,9,10}^m, \dif s_{11,12,13,14,15}^m, \dots\rangle_R\\
\langle \dif s_{1,2,3,4,5,6}^m, \dif s_{7,8,9,10,11,12}^m, \dif s_{13,14,15,16,17,18}^m, \dots\rangle_R\\
\vdots
\end{align*}
\end{lemma}
\begin{prf}
The proof proceeds in complete analogy to the proof of lemma \ref{lemma02}. By the same arguments we generalise equation \eqref{eqn02} to
\begin{align*}
\gamma(s_{i_0,j_0,k_0,\dots}^m)=\sum_{i,j} p_{i,j,k,\dots}s_{i,j,k,\dots}^0+q_{i_0,j_0,k_0,\dots}
\end{align*}
with $0\neq p_{i_0,j_0,k_0,\dots}=t_{i_0}^mt_{j_0}^mt_{k_0}^m\cdots$ (and $s_{i,j,k,\dots}^0$ of word-length $l$ as well as $q_{i_0,j_0,k_0,\dots}$ of word-length not equal to $l$) for higher word-length $l$ in the $s_i^m$. This property is proved inductively over the word-length in the $s_i^m$.

As a next step one again considers an arbitrary non-trivial linear combination
\begin{align}\label{eqn05}
v=\sum_{i} p_{i} \dif s_{i_1,i_2,i_3, i_4,\dots, i_l}^m
\end{align}
with $p_{i}\in R$ of the specified $\dif s_{i_1,i_2,i_3,i_4,\dots, i_l}^m$ and one shows that it has non-trivial image under $\gamma$. With the terminology from above we see that \linebreak[4]$\gamma(v)|_{R\otimes \Lambda^l\langle s_1^0,\dots,s_n^0\rangle}$ has $l$ ``regular terms'' and $l$ ``rest terms''. Cancellations are only possible between terms belonging to different $l$-tuples. As above we form the graph $\mathcal{G}=(V,E)$ where the vertices $V$ now correspond to the $s_{i_1,i_2,\dots, i_l}$. This graph has the very properties described in cases \eqref{case2} and \eqref{case1} (on page \pageref{case1})---in case \ref{case2-1} at most $l$ out of $2l$ terms might be cancelled, again only the ones of $x_1$---by the very same arguments. Consequently, we deduce the injectivity of $\gamma$ on the specified subspaces.
\end{prf}
We remark that the property of the formed graphs to be $3$-acyclic---i.e.~the possibility of ``localising'' the arguments---depends heavily on the fact that the integral triples belonging to the different vertices are disjoint as sets.

\begin{lemma}\label{lemma03}
Both in characteristics $2$ and $0$, the map $\gamma$ is injective on
the $R$-module
\begin{align*}
\langle
&1, s_1^m,
\\&s_{1,2}^{m}, s_{1,3}^{m}, s_{1,4}^{m}, \dots, s_{1,n}^{m},
\\&\dif s_{1,2,3}^{m}, \dif s_{4,5,6}^{m}, \dif s_{7,8,9}^{m}, \dots,
\\&s_{1,2,3}^{m}, s_{4,5,6}^{m}, s_{7,8,9}^{m}, \dots
\rangle_R
\end{align*}
\end{lemma}
\begin{prf}
As observed in \cite[Lemma 2.1.b]{Pup09}, the map $\gamma$ is injective on
\begin{align*}
V_1:=\langle 1,s_1^m,s_2^m,\dots, s_n^m, s_{1,2}^{m}, s_{1,3}^{m}, s_{1,4}^{m}, \dots, s_{1,n}^{m}\rangle_R
\end{align*}
in both characteristics $2$ and $0$.

We see that it is even injective on $V_1'\oplus V_2\oplus V_3$ where
\begin{align*}
V_1'&:=\langle 1,s_1^m,s_{1,2}^{m}, s_{1,3}^{m}, s_{1,4}^{m}, \dots, s_{1,n}^{m}\rangle_R
\intertext{and}
V_2&:=\langle\dif s_{1,2,3}^{m}, \dif s_{4,5,6}^{m}, \dif s_{7,8,9}^{m}, \dots\rangle_R
\intertext{respectively}
V_3&:=\langle s_{1,2,3}^{m}, s_{4,5,6}^{m}, s_{7,8,9}^{m}, \dots\rangle_R
\end{align*}
Indeed, in lemma \ref{lemma02} we proved its injectivity on $V_2$. Let us now see that any two elements $0\neq v_1 \in V_1'$ and $0\neq v_2\in V_2$ satisfy $\gamma(v_1)\neq \gamma(v_2)$. This, however, is obvious due to $\dif v_2=0$ and $\dif v_1\neq 0$ unless $v_1\in R$. If $v_1\in R$, the inequality $v_1\neq v_2$ is clear. This yields the injectivity on $V_1'\oplus V_2$. The injectivity on $V_1'\oplus V_2 \oplus V_3$ can be deduced as follows. Suppose that $v_1\in V_1'$, $v_2\in V_2$ and $0\neq v_3\in V_3$. We obtain $\gamma(v_1+v_2+v_3)\neq 0$, since $0\neq \dif \gamma(v_1+v_2+v_3)\in V_1'\oplus V_2$ and since $\gamma$ is injective on $V_1'\oplus V_2$, as we have seen. This yields the assertion.

\end{prf}

The last lemma paves the way towards the next one, which will provide the necessary information to improve the lower cohomological bound.
\begin{lemma}\label{lemma05}
In characteristic $0$ the map $\gamma$ is injective on
the $R$-module
\begin{align*}
\langle
&1,
\\&s_1^m,s_2^m,\dots, s_n^m,
\\&s_{1,2}^{m}, s_{1,3}^{m}, s_{1,4}^{m}, \dots, s_{1,n}^{m},
\\&\dif s_{1,2,3}^{m}, \dif s_{4,5,6}^{m}, \dif s_{7,8,9}^{m}, \dots,
\\&s_{1,2,3}^{m}, s_{4,5,6}^{m}, s_{7,8,9}^{m}, \dots
\rangle_R
\end{align*}
if $\gamma$ is degree-preserving.
\end{lemma}
\begin{prf}
We use the notation from the proof of lemma \ref{lemma03}. In characteristic $0$ we even obtain the injectivity on $V_1\oplus V_2\oplus V_3$
by a similar but slightly extended reasoning. Indeed, for non-trivial $v_1\in V_1$, $v_2\in V_2$ and $v_3\in V_3$ we proceed as follows.
In characteristic zero we have $\deg t_i=2$ and, consequently, $\deg s_i^m=2m+1$. Thus we may split $v_1=v_1'+v_1''$ with $v_1'\in \langle s_1^m,s_2^m,\dots, s_n^m\rangle_R$ (consisting of summands) of odd degree, $v_1''\in  \langle 1,s_{1,2}^{m}, s_{1,3}^{2m}, s_{1,4}^{m}, \dots, s_{1,n}^{m}\rangle_R$ (consisting of summands) of even degree.
That is, we have a $\zz_2$-grading on $K_i(m)$. In the assertion we assume the map $\gamma$ to preserve this grading, in particular. This implies that, given the injectivity on the odd part, i.e.~$\gamma(v_1'+v_3)\neq 0$, and the injectivity on the even part, i.e.~$\gamma(v_1''+v_2)\neq 0$, we derive the injectivity on $V_1\oplus V_2\oplus V_3$. The fact that $\gamma(v_1''+v_2)\neq 0$ follows as in the proof of lemma \ref{lemma03}.

The inequality $\gamma(v_1'+v_3)\neq 0$ is due to the following arguments adapted from \cite[Lemma 2.3]{Pup09}. We may assume $v_1'+v_3$ to be a homogeneous element. We write
\begin{align*}
v_1'+v_3=\sum_i p_i s_{i_1,i_2,i_3}^{m}+\sum_j p_j s_{j}^m
\end{align*}
with $p_i, p_j\in R$. In order to show that $\gamma(v_1'+v_3)\neq 0$ we form the ring $R/I$ where $I$ is the ideal generated by the $t_{i_1}^{m+1} p_{i}, t_{i_2}^{m+1} p_{i}, t_{i_3}^{m+1} p_{i}, t_j^{m+1} p_j$. In the cohomology with coefficients in $R/I$ the element $v_1'+v_3$ becomes a cycle. As observed in \cite{Pup09} the map $\gamma$ and the \emph{multiplicative} morphism $\iota\co K_n(m)\to K_n(0)$ induced by $\iota(s_j^m)=t^ms_j^0$ are homotopic. Thus, in particular, we obtain $H(\gamma,R/I)=H(\iota,R/I)$. We shall now see that $H(\iota, R/I)([v_1'+v_3])\neq 0$, which will complete the proof.

For this we compute
\begin{align*}
\iota(v_1'+v_3)&=\iota\bigg(\sum_i p_i s_{i_1,i_2,i_3}^{m}+\sum_j p_j s_{j}^m\bigg)\\&=\sum_i p_i t_{i_1}^mt_{i_2}^mt_{i_3}^m s_{i_1}^0s_{i_2}^0s_{i_3}^0 + \sum_j p_j t_j^m s_j^0
\end{align*}
We have that $\deg t_{i_1}^mt_{i_2}^mt_{i_3}^m s_{i_1}^0s_{i_2}^0s_{i_3}^0=6m+3$ and $\deg t_j^m s_j^0=2m+1$, which implies that $\deg p_j=\deg p_i+4m+2$. Consequently,
\begin{align*}
\deg p_i s_{i_1,i_2,i_3}^{m}&=(\deg p_j-4m-2)+(6m+3)\\&=\deg p_j+2m+1\\&<2(m+1)\deg p_j\\&=\deg t_j^{m+1}p_j
\end{align*}
and the $p_i s_{i_1,i_2,i_3}^{m}$ are not in the ideal generated by the $t_j^{m+1} p_j$.

In other words, this proves that $\gamma(v_1')=-\gamma(v_3)$ implies $[\iota(v_1')]=-[\iota(v_3)]$, which itself yields that $H(\iota,R/J)([v_3])=0$ where $J$ is the ideal generated in $R$ by the $t_{i_1}^{m+1} p_{i}, t_{i_2}^{m+1} p_{i}, t_{i_3}^{m+1} p_{i}$. This implies that $H(\gamma,R/J)([v_3])=0$, since $\gamma$ and $\iota$ are homotopic. Finally, we deduce $\gamma(v_3)=0$, since the ideal $J$ is generated in degrees larger than $\deg v_3$. By the injectivity of $\gamma$ on $V_3$ we infer that $v_3=0$. Combining this with the injectivity on $V_1$ finishes the proof.
\end{prf}
\begin{rem}
One might be tempted to generalise the cohomological reasoning at the end of the proof to showing that $\gamma$ is injective on more levels than just $1$ and $3$. That is, one would like to show that it is injective on the direct sum of subspaces in several word-lenghts in the $s_i^m$---not necessarily only $1$ and $3$---if and only if so it is on each respective subspace. However, the next levels are $4m+2$, $8m^2+6m+1$, etc. In the main application $m$ will correspond to the dimension of the orbit space of an almost free $\s^1$-torus action. Since the maximal rank of an $\s^1$-torus acting freely on an $m$-dimensional finite CW-complex is $m$, these considerations will be of no use.
\end{rem}
\begin{rem}
Unfortunately, it is too much to hope for the injectivity on the specified space also in characteristic $2$. In \cite[Example 2.2]{Pup09} a morphism $\gamma\co K_3(1)\to K_3(0)$ is provided which is not injective on $\langle \dif s_{1,2,3}^1, \dif s_{1,2}^1\rangle_R$.
\end{rem}
Recall that the rank $\rk \gamma$ of $\gamma$ is defined to be the rank of the induced $R_{(0)}$-linear morphism when localising $R$ to its field of fractions $R_{(0)}$.
\begin{cor}\label{cor01}
In characteristic $0$ we have that $\rk \gamma\geq 2(n+\lfloor n/3\rfloor)$ for a degree-preserving $\gamma$.
\end{cor}
\begin{prf}
Since localisation is an exact functor, the injectivity of $\gamma$ on
\begin{align*}
\langle &1,s_1^m,s_2^m,\dots, s_n^m, s_{1,2}^{m}, s_{1,3}^{m}, s_{1,4}^{m}, \dots, s_{1,n}^{m},\\& \dif s_{1,2,3}^{m}, \dif s_{4,5,6}^{m}, \dots, s_{1,2,3}^{m}, s_{4,5,6}^{m}, \dots\rangle_R
\end{align*}
from lemma \ref{lemma05} implies the injectivity of $\gamma\otimes R_{(0)}$ on the $R_{(0)}$-vector space
\begin{align*}
\langle &1,s_1^m,s_2^m,\dots, s_n^m, s_{1,2}^{m}, s_{1,3}^{m}, s_{1,4}^{m}, \dots, s_{1,n}^{m}, \\&\dif s_{1,2,3}^{m}, \dif s_{4,5,6}^{m},\dots, s_{1,2,3}^{m}, s_{4,5,6}^{m}, \dots\rangle_{R_{(0)}}
\end{align*}
Since $\dif s_{i_1,i_2,i_3}^{m}=t_{i_1}^{m+1}s_{i_2,i_3}^{m}\pm t_{i_2}^{m+1}s_{i_1,i_3}^{m}\pm t_{i_3}^{m+1}s_{i_1,i_2}^{m}$, we can count as follows: We have the $2n$ basis elements $1,s_1^m,s_2^m,\dots, s_n^m, s_{1,2}^{m}, s_{1,3}^{m}, s_{1,4}^{m}, \dots, s_{1,n}^{m}$. Since each of the $\dif s_{i_1,i_2,i_3}^m$ contributes (at least one) new basis element $s_{i_2,i_3}$, this adds $\lfloor n/3\rfloor$ to the dimension. Finally, taking into account the $\lfloor n/3\rfloor$ additional elements $s_{i_1,i_2,i_3}$, we see that
the space on which the localisation of $\gamma$ is injective has dimension $2(n+\lfloor n/3\rfloor)$ over $R_{(0)}$. This yields the result.
\end{prf}

\section{Proof of theorem A} \label{sec02}

In this section---following \cite{Pup09}---we shall present the topological background which will permit us to derive theorem \ref{theoA}.

For this it is necessary to recall the \emph{minimal Hirsch--Brown model}---see \cite[sections 3,5, p.~372]{Pup09}, \cite[section 1.3, p.~22]{AP93}. It is given by
\begin{align*}
H^*(X,\Bbbk)\tilde \otimes H^*(\B G,\Bbbk)
\end{align*}
where the tensor product carries a twisted differential and where the group $G$ is either an $\s^1$-torus or a $\zz_2$-torus. From now on let us focus on the case of $\operatorname{char} \Bbbk=0$. This model algebraically captures the information given in the Borel construction and can be compared to the Leray--Serre spectral sequence of the fibration $X\hto{}X_G\to \B G$ respectively its corresponding Sullivan model in Rational Homotopy Theory. The deviation of the tensor product from being non-twisted reflects the deviation of the fibration from being a Cartesian product.

We quote the following famous
\begin{theo}[Hsiang]\label{theo01}
Let $G$ be a compact Lie group acting continuously on $X$. Then $G$ acts almost freely if and only if $H^*_G(X,\qq)$ is finite-dimensional.
\end{theo}
(see~\cite{Hsi75}, \cite{AH78}). If $G$ acts almost freely, we obtain a weak homotopy equivalence $X_G\simeq X/G$ and
\begin{align}\label{eqn06}
H^*(X_G,\Bbbk)\cong H(H^*(X,\Bbbk)\tilde \otimes H^*(\B G,\Bbbk))\cong H^*(X/G,\Bbbk)
\end{align}
The two main technical tools to connect the Hirsch--Brown model to the Koszul complexes are provided by the subsequent lemmas, which we cite from \cite[Proposition 4.1, Proposition 4.2]{Pup09} respectively from \cite[Lemma 1.4.17, Lemma 1.4.18]{AP93} in the graded context.

We consider free cochain complexes $\tilde C$ over $R$ with an augmentation $\varepsilon\co \tilde C\to \Bbbk$ which induces a surjection in cohomology.
\begin{lemma}\label{lemma07}
Let $\tilde C$ be as above and assume that
$H^{>m}(\tilde C)$ vanishes. Then there exists a map of $R$-complexes $\alpha\co K_n(m) \to \tilde C$
such that
\begin{align*}
H^*(\varepsilon)H^*(\alpha)\co H(K_n(m))=R/(t_1^{m+1},\dots, t_r^{m+1})\to \Bbbk
\end{align*}
is the canonical projection.
\end{lemma}
Suppose now that $\tilde C$ admits a filtration $\mathcal{F}_*(\tilde C)$ which satisfies the following:
\begin{itemize}
\item There are proper inclusions of free subcomplexes
\begin{align*}
\mathcal{F}_0(\tilde C)=0 \In \mathcal{F}_1(\tilde C)\In  \ldots \In \mathcal{F}_{\ell(\mathcal{F}_*(\tilde C))}(\tilde C)=\tilde C
\end{align*}
which are direct summands as $R$-modules and
\item it holds that $\bar \dif(\mathcal{F}_i(\tilde C))\In \mathcal{F}_{i-1}(\tilde C)$ for $i=0, \dots, \ell(\mathcal{F_*(\tilde C)})$.
\end{itemize}
where $\ell(\mathcal{F_*(\tilde H)})$ is the length of the filtration.
\begin{lemma}\label{lemma08}
Then there exists a filtration-preserving morphism of \linebreak[4]$R$-complexes $\beta\co \tilde C\to K_n(0)$ which commutes with the respective augmentations.
\end{lemma}

We shall now apply these lemmas to the Hirsch--Brown model from above. Due to theorem \ref{theo01} and the isomorphism \eqref{eqn06}, we obtain that the cohomology of the Hirsch--Brown model vanishes above the formal dimension $m:=\dim X/G<\infty$ of the orbit space. By the general theory, we then obtain a unital morphism
\begin{align*}
\beta\circ\alpha\co K_n(m) \to K_n(0)
\end{align*}
factoring through the minimal Hirsch--Brown model.

\begin{proof}[\textsc{Proof of theorem \ref{theoA}}]
Since the map $\beta\circ \alpha$ factors through the minimal Hirsch--Brown model, we deduce that
\begin{align*}
\dim_{\Bbbk} H^*(X,\Bbbk)=\operatorname{rk}_R (H^*(X,\Bbbk)\tilde \otimes H^*(\B G,\Bbbk))\geq \rk \gamma
\end{align*}
Corollary \ref{cor01} shows that $\rk \gamma\geq 2(n+\lfloor n/3\rfloor)$, which yields the result.

\vspace{5mm}

Alternatively, we might use the injectivity of $\gamma$ on the $(n+\lfloor n/3\rfloor)$-dimen\-sional space
\begin{align*}
\langle &1, s_{1,2}^{m}, s_{1,3}^{m}, s_{1,4}^{m}, \dots, s_{1,n}^{m},\dif s_{1,2,3}^{m}, \dif s_{4,5,6}^{m}, \dots\rangle_{R_{(0)}}
\end{align*}
to argue as follows: The cohomology of the localised Hirsch--Brown model, i.e.~when passing from the $R$-module to the $R_{(0)}$-vector space, vanishes, since
$1$ becomes a boundary. When localising, the $\zz$-grading induces a $\zz_2$-grading into even and odd part. We derive that the dimension of the odd-degree part equals the dimension of the even-dimensional part. Since the space above is concentrated in even degrees, the sum over all even-degree Betti numbers of $H^*(X,\Bbbk)$ must be at least $n+\lfloor n/3\rfloor$, since $\beta\circ\alpha$ factors through the minimal Hirsch--Brown model. Thus the dimension of $H^*(X,\Bbbk)$ is at least $2(n+\lfloor n/3\rfloor)$.
\end{proof}



\vfill

\begin{center}
\noindent
\begin{minipage}{\linewidth}
\small \noindent \textsc
{Manuel Amann} \\
\textsc{Department of Mathematics}\\
\textsc{University of Pennsylvania}\\
\textsc{David Rittenhouse Laboratory}\\
\textsc{209 South 33rd Street}\\
\textsc{Philadelphia, Pennsylvania}\\
\textsc{PA 19104-6395} \\
\textsc{USA}\\
[1ex]
\textsf{mamann@uni-muenster.de}\\
\textsf{http://hans.math.upenn.edu/$\sim$mamann/}
\end{minipage}
\end{center}

\end{document}